\documentclass[12pt]{article}
\usepackage{amsmath,amssymb}
\def\L{\mathcal{L}}
\def\A{\mathbf{A}}
\def\B{\mathbf{B}}
\def\a{\mathbf{a}}
\def\b{\mathbf{b}}
\newtheorem{theorem}{\hspace*{\parindent}Theorem}
\newtheorem{lemma}{\hspace*{\parindent}Lemma}
\newtheorem{corollary}{\hspace*{\parindent}Corollary}
\newtheorem{conjecture}{\hspace*{\parindent}Conjecture}
\newcounter{theremark}
\newcommand{\rem}{\par\refstepcounter{theremark}\textbf{Remark \arabic{theremark}.} }
\title{Completely monotonic gamma ratio and infinitely divisible H-function of Fox}
\author{D.B.\:Karp$^{\rm a,c}$\footnote{Corresponding author. E-mail: D. Karp -- \emph{dimkrp@gmail.com}, E.\:Prilepkina --  \emph{pril-elena@yandex.ru}}~~and E.G.\:Prilepkina$^{\rm a,b}$
\\[10pt]\small{\textit{$\phantom{1}^a$Far Eastern Federal University, 8 Sukhanova street, Vladivostok, 690950, Russia}}\\\small{\textit{$\phantom{1}^b$Institute of Applied Mathematics, FEBRAS, 7 Radio Street, Vladivostok,  690041, Russia}}
\\\small{\textit{$\phantom{1}^c$Universidad del Atl\'{a}ntico, Km~7 Antigua via, Puerto Colombia, Colombia}}}
\date{}
\begin{document}
\maketitle

\begin{center}
\parbox{12cm}{
\small\textbf{Abstract.}
We investigate conditions for logarithmic complete monotonicity of a quotient of two products of gamma functions, where the argument of each gamma function has different scaling factor.  We give necessary and sufficient conditions in terms of nonnegativity of some elementary function and more practical sufficient conditions in terms of parameters.  Further, we study the representing measure in Bernstein's theorem for both equal and non-equal scaling factors.  This leads to conditions on parameters under which Meijer's $G$-function or Fox's $H$-function represents an infinitely divisible probability distribution on positive half-line.  Moreover, we present new integral equations for both $G$-function and $H$-function. The results of the paper generalize those due to Ismail (with Bustoz, Muldoon and Grinshpan) and Alzer who considered previously the case of unit scaling factors.
}
\end{center}

\bigskip


\bigskip


\bigskip
\bigskip

\paragraph{1. Introduction.}
Recall that a nonnegative function $f$ defined on $(0,\infty)$ is
called completely monotonic ({\bf c.m.}) if it has derivatives of
all orders and $(-1)^nf^{(n)}(x)\ge0$ for $n\geq 1$ and $x>0$
\cite[Definition~1.3]{SSV}. This inequality  is known to be strict
unless $f$ is a constant. By the celebrated Bernstein theorem a
function is completely monotonic if and only if it is the Laplace
transform of a nonnegative measure \cite[Theorem~1.4]{SSV}.
 The above definition implies the following equivalences
\begin{multline}\label{eq:cm-equiv}
f  \text{ is c.m. on } (0,\infty) \Leftrightarrow f\geq 0 \text{ and } -\!f' \text{ is c.m. on } (0,\infty)
\\
\Leftrightarrow -f'\text { is c.m. on } (0,\infty) \text{ and } \lim\limits_{x\to\infty} f(x)\geq 0.
\end{multline}

A positive function $f$ is said to be logarithmically completely
monotonic ({\bf l.c.m.}) if $-(\log{f})'$ is completely monotonic
\cite[Definition~5.8]{SSV}. According to (\ref{eq:cm-equiv})
\begin{multline}\label{eq:lcm}
f \text{ is l.c.m. on } (0,\infty) \Leftrightarrow (-\log f(x))'\ge0 \text{ and } (\log{f})'' \text{ is c.m. on } (0,\infty)
\\
\Leftrightarrow (\log f)'' \text{ is c.m. on } (0,\infty) \text{ and } \lim\limits_{x\to\infty} (-\log f(x))'\geq 0.
\end{multline}
The class of l.c.m. functions is a proper subset of the class of
c.m. functions. Their importance stems from the fact that they
represent Laplace transforms of infinitely divisible probability
distributions, see \cite[Theorem~5.9]{SSV} and
\cite[Section~51]{Sato}.

 The study
of complete monotonicity of the ratio
\begin{equation}\label{eq:classical-ratio}
U(x)=\prod_{i=1}^{p}\frac{\Gamma(x+a_i)}{\Gamma(x+b_i)},
\end{equation}
where $\Gamma$ stands for Euler's gamma function and $p\geq 1$,
has been initiated by Bustoz and Ismail  who demonstrated in their
1986 paper \cite{BusIsm} that for $p=2$, $a_1=0$ and
$a_1+a_2=b_1+b_2$ this function is l.c.m. on $(0,\infty)$. Eight
years later Ismail and Muldoon showed in \cite{IsmailMuldoon94}
that $U(x)$ is l.c.m. on $(0,\infty)$ for general $p$ if
$\sum\nolimits_{i=1}^p{a_i}=\sum\nolimits_{i=1}^p {b_i}$ and
$b_1=b_2=\cdots=b_p>0$, $a_i\ge0$, or $b_1=b_2=\cdots=b_{p-1}=0$,
$b_p>0$, $a_i\ge0$. Their result is, in fact, formulated for the
ratio of $q$-gamma functions, but the proof works for $U(x)$ just
as well.  The subject was further pursued by Alzer \cite{Alzer}  who  showed in
1997 that $U(x)$ is l.c.m. on $(0,\infty)$ if
\begin{equation}\label{eq:amajorb}
\begin{split}
& 0\leq{a_1}\leq{a_2}\leq\cdots\leq{a_p},~~
0\leq{b_1}\leq{b_2}\leq\cdots\leq{b_p},
\\
&~\text{and}~\sum\limits_{i=1}^{k}a_i\leq\sum\limits_{i=1}^{k}b_i~~\text{for}~~k=1,2\ldots,p.
\end{split}
\end{equation}
These inequalities are known as weak supermajorization \cite[Definition~A.2]{MOA} and are abbreviated as
$\mathbf{b}\!\prec^W\!\!\mathbf{a}$, where $\mathbf{a}\!=\!(a_1,\ldots,a_p)$,
$\mathbf{b}\!=\!(b_1,\ldots,b_p)$. In their 2006 paper \cite{GrinIsm} Grinshpan and Ismail found new sufficient
conditions for $U(x)$ to be l.c.m. of  when $p=2^n$ or $p=n!/2$, $n\geq 1$.  We will explain and slightly
generalize their results in the next section.  Finally, in 2009 Guo and Qi \cite{Guo-Qi} used another approach to investigate
logarithmic complete monotonicity of $U(x)$ for arbitrary real values of $a_i$, $b_i$.  Their results, however, lead to complete
monotonicity of $U(x)$ on some subinterval of $(0,\infty)$ of the form $(\gamma,\infty)$, where $\gamma>0$.

The purpose of the present  paper is to investigate complete monotonicity of the ratio
\begin{equation}\label{eq:Gammaratio}
W(x)=\frac{\prod_{i=1}^{p}\Gamma(A_ix+a_i)}{\prod_{j=1}^{q}\Gamma(B_jx+b_j)},
\end{equation}
where $\A=(A_1,\ldots,A_p)$ and $\B=(B_1,\ldots,B_q)$ are strictly positive scaling factors, while $\a=(a_1,\ldots,a_p)$ and  $\b=(b_1,\ldots,b_q)$ are nonnegative.  We find a necessary and sufficient condition for $W$ to be l.c.m. on $(0,\infty)$ in terms of nonnegativity of some function which, unfortunately, is not easy to verify.  Further, we supply several simple sufficient conditions for such nonnegativity in terms of the vectors $\A,\B,\a,\b$ as well as some necessary conditions.  When $W$ is completely monotonic, we proceed by deriving its representing measure in the Bernstein theorem which leads us to studying some new properties of Fox's $H$-function.  We begin, however, by revisiting the ratio defined in (\ref{eq:classical-ratio}) which we call the unweighted case. In the following section we refine some of the known results for $U(x)$ and discuss its representing measure in Bernstein's theorem.

\bigskip

\paragraph{2. The unweighted case revisited.}
All proofs of complete monotonicity of $U(x)$ we are aware of are
based on the integral representation of  the logarithmic
derivative of the gamma function \cite[p.16]{MOS}
\begin{equation}\label{eq:psi-Gauss}
\psi(z)=\frac{\Gamma'(z)}{\Gamma(z)}=-\gamma+\int_{0}^{\infty}\frac{e^{-t}-e^{-tz}}{1-e^{-t}}dt,
\end{equation}
where $\gamma$ is Euler-Mascheroni constant.  The next two properties of digamma function $\psi$ will be repeatedly used in the sequel
\begin{equation}\label{eq:psi-deriv}
\psi'(x)=\int\limits_{0}^{\infty}\frac{te^{-xt}}{1-e^{-t}}dt,~~x>0,
\end{equation}
and
\begin{equation}\label{eq:psi-asimptotic}
 \psi(z)=\log(z)-\frac{1}{2z}+O(z^{-2})~~\text{as}~~z\to\infty,
\end{equation}
see \cite[Theorem~1.6.1,  Corollary~1.4.5]{AAR}.

Representation (\ref{eq:psi-Gauss}) lead Grinshpan and
Ismail to formulate Lemma~2.1 in \cite{GrinIsm}  stating that
$q(x)=\prod_{k=1}^{p}\Gamma^{\sigma_k}(x+\lambda_k)$ with
$\sum_{k=1}^{p}\sigma_k=0$ and $\lambda_k\ge0$ is l.c.m. on
$(0,\infty)$ if and only if
$$
v(t)=\sum_{k=1}^{p}\sigma_kt^{\lambda_k}\ge0~~\text{for}~t\in(0,1].
$$
In fact, looking at asymptotic expansion of $-(\log{q(x)})'$ as
$x\to\infty$ it is easy to sharpen this  lemma as follows: $q(x)$
is l.c.m. on $(0,\infty)$ iff $v(t)\ge0$ on $(0,1]$ and
$\sum_{k=1}^{p}\sigma_k=0$.  When $q(x)=U(x)$, where $U$ is defined in (\ref{eq:classical-ratio}), we have
\begin{equation}\label{eq:munts}
U(x)\ \mathrm{is \ l.c.m.\ on}\ (0,\infty)\Leftrightarrow
v(t)=\sum_{k=1}^{p}(t^{a_k}-t^{b_k})\geq 0\ \mathrm{on}\ (0,1].
\end{equation}
As mentioned above Alzer noticed that $v(t)\ge0$ on $(0,1]$ under majorization conditions (\ref{eq:amajorb}) which follows from 1949 result
of Tomi\'{c} \cite[Proposition~4.B.2]{MOA}. Let us remark that for $p=2$ conditions $0\leq\min(a_1,a_2)\leq\min(b_1,b_2)$ and $a_1+a_2\le{b_1+b_2}$,
equivalent to (\ref{eq:amajorb}), are necessary and sufficient for nonnegativity of $v(t)$.
 Necessity follows by considering the asymptotics of $v(t)$ as $t\to0$ and $t\to1$.
A proof of a more general result is given in Corollary~\ref{cr:necessary} below.

Grinshpan and Ismail took another path in \cite{GrinIsm} and considered two types of factorization
$$
v(t)=\prod_{1\le{i}<j\le{n}}(t^{\alpha_j}-t^{\alpha_i})~~\text{and}~~v(t)=\prod_{i=1}^{n}(1-t^{\alpha_i}).
$$
It clear that for $\alpha_1>\alpha_2>\cdots>\alpha_n>0$ for the first factorization and $\alpha_i>0$ for the second
we get expressions nonnegative on $(0,1]$.  The corresponding value of $p$ is $n!/2$ for the first factorization and $2^{n-1}$ for the second.  In both cases the authors of \cite{GrinIsm}  found explicit combinatorial descriptions of the vectors $\a$, $\b$ in (\ref{eq:classical-ratio}) which lead to the above factorizations.  We would like to remark here that taking $v$ in the form
$$
v(t)=\prod_{i=1}^{n}(t^{\beta_i}-t^{\alpha_i})
$$
with $\alpha_i\ge\beta_i\ge0$ for $i=1,\ldots,n$ certainly results in nonnegative $v$ on $(0,1]$. This, of course, can be reduced to
the second factorization of Grinshpan and Ismail  by factoring out $t^{\sum\beta_i}$. Also one can get the first type of factorization by properly choosing $\beta_i$.  The vector $\a$ ($\b$) is expressed in \cite{GrinIsm}  in terms of sum over even (odd) permutations of the numbers $\alpha_i$ for the first factorization, and in terms of sums of the numbers $\alpha_i$ indexed by components of increasing integer vectors for the second factorization.  In the theorem below which generalizes both factorizations, the description of the vectors $\a$ and $\b$ in terms of the numbers $\alpha_i$ and $\beta_i$ turns out to be different and simpler.
\begin{theorem}\label{th:unweighted}
Let $I=\{1,2,\ldots,n\}$ and suppose $\mathcal{I}_\text{odd}$ \emph{(}\,$\mathcal{I}_{\text{even}}$\,\emph{)}
comprises all subsets of $I$ having odd \emph{(}even\emph{)} number of elements, $\emptyset\in\mathcal{I}_{\text{even}}$.
Suppose that $\alpha_i\ge\beta_i\ge0$ for $i=1,\ldots,n$.  Then
$$
U(x)=\frac{\prod\nolimits_{J\in\mathcal{I}_{\text{even}}}\Gamma\left(x+\sum_{i\in{J}}\alpha_i+\sum_{i\in{I_n\setminus{J}}}\beta_i\right)}
{\prod\nolimits_{J\in\mathcal{I}_{\text{odd}}}\Gamma\left(x+\sum_{i\in{J}}\alpha_i+\sum_{i\in{I_n\setminus{J}}}\beta_i\right)}
$$
is l.c.m. on $(0,\infty)$.
\end{theorem}

\textbf{Proof.}  We will write $I=I_n$, $\mathcal{I}_\textit{odd}=\mathcal{I}_\textit{odd}(n)$, $\mathcal{I}_{\textit{even}}=\mathcal{I}_{\textit{even}}(n)$ to emphasize the dependence on $n$. According to (\ref{eq:munts}) it suffices to verify that  $v(t)\ge0$ on $(0,1]$, where
$$
v(t)=\sum\limits_{J\in\mathcal{I}_{\textit{even}}(n)}t^{\sum\limits_{i\in{J}}\alpha_i+\sum\limits_{i\in{I_n\setminus{J}}}\beta_i}-
\sum\limits_{J\in\mathcal{I}_{\textit{odd}}(n)}t^{\sum\limits_{i\in{J}}\alpha_i+\sum\limits_{i\in{I_n\setminus{J}}}\beta_i}.
$$
We will demonstrate that
$$
v(t)=\prod_{i=1}^{n}(t^{\beta_i}-t^{\alpha_i})
$$
by induction in $n$. Indeed, the case $n=1$ is obvious.  The induction step now follows from the identity
$$
\sum\limits_{J\in\mathcal{I}_{\textit{even}(n)}}t^{\sum\limits_{i\in{J}}\alpha_i+\sum\limits_{i\in{I_n\setminus{J}}}\beta_i}=
\sum\limits_{J\in\mathcal{I}_{\textit{even}}(n-1)}t^{\sum\limits_{i\in{J}}\alpha_i+\sum\limits_{i\in{I_{n-1}\setminus{J}}}\beta_i+\beta_n}
+\sum\limits_{J\in\mathcal{I}_{\textit{odd}}(n-1)}t^{\sum\limits_{i\in{J}}\alpha_i+\alpha_n+\sum\limits_{i\in{I_{n-1}\setminus
J}}\beta_i},
$$
and similar identity for the negative term in $v(t)$.  Nonnegativity is now obvious from $\alpha_i\ge\beta_i\ge0$ for $i=1,\ldots,n$. $~\hfill \square$

The next example generalizes \cite[Corollary~1.4]{GrinIsm}

\textbf{Example~1.} For $\alpha_i\ge\beta_i\ge0$ for $i=1,2,3$ the function
$$
U(x)=\frac{\Gamma(x+\beta_1+\beta_2+\beta_3)\Gamma(x+\beta_1+\alpha_2+\alpha_3)\Gamma(x+\alpha_1+\beta_2+\alpha_3)\Gamma(x+\alpha_1+\alpha_2+\beta_3)}
{\Gamma(x+\alpha_1+\beta_2+\beta_3)\Gamma(x+\beta_1+\alpha_2+\beta_3)\Gamma(x+\beta_1+\beta_2+\alpha_3)\Gamma(x+\alpha_1+\alpha_2+\alpha_3)}
$$
is l.c.m. on $(0,\infty)$.

Further, changing variable in \cite[formula~2.24.2.1]{PBM3} or in particular case of \cite[Theorem~2.2]{KilSaig} we get
\begin{equation}\label{eq:G-represent}
\prod\limits_{i=1}^{p}\frac{\Gamma(x+a_i)}{\Gamma(x+b_i)}=\int\limits_{(0,\infty)}e^{-tx}
G^{p,0}_{p,p}\left(e^{-t}\,\,\vline\begin{array}{c}\!\b\\
\!\a\end{array}\!\!\right)\!dt
\end{equation}
provided that $\sum_{i=1}^{p}(b_i-a_i)>0$.  Here $G^{p,0}_{p,p}$
is Meijer's $G$-function defined by the contour integral
\begin{equation*}\label{eq:G-defined}
G^{p,0}_{p,p}\!\left(\!z~\vline\begin{array}{l}\b\\\a\end{array}\!\!\right)
=\frac{1}{2\pi{i}}
\int\limits_{\L}\!\!\frac{\Gamma(a_1\!+\!s)\dots\Gamma(a_p\!+\!s)}
{\Gamma(b_{1}\!+\!s)\dots\Gamma(b_p\!+\!s)}z^{-s}ds.
\end{equation*}
This function is a particular case of Fox's $H$-function defined in (\ref{eq:Fox}) below on setting $p=q$ and $A_i=B_i=1$,
$i=1,\ldots,p$ in that definition.  See text below (\ref{eq:Fox}) for a description of the contour $\L$.
It follows from an expansion due to N{\o}rlund \cite[(2.28)]{Norlund} that for $\sum_{i=1}^{p}(b_i-a_i)=0$ formula (\ref{eq:G-represent}) must be modified as
\begin{equation}\label{eq:psi0measure}
\prod\limits_{i=1}^{p}\frac{\Gamma(x+a_i)}{\Gamma(x+b_i)}=\int\limits_{[0,\infty)}e^{-tx}
\left\{\delta_0+G^{p,0}_{p,p}\left(e^{-t}\,\,\vline\begin{array}{c}\!\b\\
\!\a\end{array}\!\!\right)\right\}\!dt,
\end{equation}
where $\delta_0$ denotes the unit mass  concentrated at zero.  More information about N{\o}rlund's expansions can be found in our forthcoming paper \cite{KPNew}.  Representation (\ref{eq:G-represent}) was previously observed by us in \cite{KPJMAA12}.

The l.c.m $\Rightarrow$ c.m. implication mentioned above $U$ translates into the next assertion:
$$
v(t)=\!\!\sum_{k=1}^{p}(t^{a_k}-t^{b_k})\ge0~\text{on}~[0,1)~\Rightarrow~G^{p,0}_{p,p}\left(x\,\,\vline\begin{array}{c}\!\b\\\!\a\end{array}\!\!\right)\ge0~\text{on}~(0,1)
$$
under additional conditions $a_i\ge0$, $i=1,\ldots,p$, and $\sum_{i=1}^{p}(b_i-a_i)>0$. The following much stronger assertion
is supported by numerical evidence:

\begin{conjecture}\label{cj:zerosG}
If $a_i\ge0$, $i=1,\ldots,p$ and $\sum_{i=1}^{p}(b_i-a_i)>0$ then
$$
\#\biggl\{\text{zeros of}~G^{p,0}_{p,p}\left(x\,\,\vline\begin{array}{c}\!\b\\\!\a\end{array}\!\!\right)~\text{on}~(0,1)\biggr\}
\le \#\biggl\{\text{zeros of}~v(t)~\text{on}~(0,1)\biggr\}.
$$
\end{conjecture}
We think that this conjecture is encoded in the following integral equation for Meijer's $G$-function which we believe to be new.
\begin{theorem}\label{th:G-integral}
Suppose $\a,\b\ge0$ and $\psi>0$. Then Meijer's $G$-function satisfies the following integral equation\emph{:}
$$
\log(1/x)G^{p,0}_{p,p}\left(x\,\,\vline\begin{array}{c}\!\b\\\!\a\end{array}\!\!\right)
=\int_x^{1}G^{p,0}_{p,p}\left(t\,\,\vline\begin{array}{c}\!\b\\\!\a\end{array}\!\!\right)
\sum_{k=1}^{p}\left(\frac{x^{a_k}}{t^{a_k}}-\frac{x^{b_k}}{t^{b_k}}\right)\frac{dt}{t-x}
$$
for $0<x<1$.
\end{theorem}
We postpone the proof of this result until proof of Theorem~\ref{th:hrepresentation}, containing a more general statement.
Using (\ref{th:G-integral}) we were also able to demonstrate that if $G^{p,0}_{p,p}(x)$ has a zero on $(0,1)$ then $v(t)$ also has at least one zero on this interval.

\bigskip

\paragraph{3. Generalized gamma ratio.}

For the positive vectors $\A=(A_1,\ldots,A_p)$,
$\B=(B_1,\ldots,B_q)$ and nonnegative vectors $\a=(a_1,\ldots,a_p)$, $\b=(b_1,\ldots,b_q)$,
consider the positive function
$$
W(x)=\frac{\prod_{i=1}^{p}\Gamma(A_ix+a_i)}{\prod_{j=1}^{q}\Gamma(B_jx+b_j)}.
$$
defined  on $(0,\infty)$.  Our goal is to study logarithmic complete monotonicity of $W$.  It is convenient to start with
\begin{lemma}\label{lm:lcm-derivative}
The function $(\log W)''$ is c.m. if and only if
\begin{equation}\label{eq:Pu-def}
P(u)=\sum\limits_{i=1}^{p}\frac{e^{-a_iu/A_i}}{1-e^{-u/A_i}}-\sum\limits_{i=1}^{q}\frac{e^{-b_iu/B_i}}{1-e^{-u/B_i}}\ge0~~\text{for
all}~u>0.
\end{equation}
In the affirmative case
\begin{equation}\label{eq:represent 1}
(\log{W})''=\int\limits_{0}^{\infty}e^{-xu}uP(u)du.
\end{equation}
\end{lemma}

\textbf{Proof.} Differentiating (\ref{eq:psi-Gauss}) we have
$$
(\log{W})''=\sum\limits_{i=1}^{p}A_i^2\psi'(A_ix+a_i)-\sum\limits_{j=1}^{q}B_j^2\psi'(B_jx+b_j).
$$
Using (\ref{eq:psi-deriv}) and making change of variable $u=A_it$ or $u=B_it$ in the appropriate integrals we get
 \begin{multline*}
(\log{W})'' =
\sum\limits_{i=1}^{p}A_i^2\int_{0}^{\infty}\frac{te^{-(A_ix+a_i)t}dt}{1-e^{-t}}
-\sum\limits_{j=1}^{q}B_j^2\int_{0}^{\infty}\frac{te^{-(B_jx+b_j)t}dt}{1-e^{-t}}
\\
=\sum\limits_{i=1}^{p}\int_{0}^{\infty}\frac{ue^{-xu-ua_i/A_i}du}{1-e^{-u/A_i}}-
\sum\limits_{j=1}^{q}\int_{0}^{\infty}\frac{ue^{-xu-ub_j/B_j}du}{1-e^{-u/B_j}}
\\
=
\int_{0}^{\infty}e^{-xu}u\left\{\sum\limits_{i=1}^{p}\frac{e^{-a_iu/A_i}}{1-e^{-u/A_i}}
-\sum\limits_{j=1}^{q}\frac{e^{-b_ju/B_j}}{1-e^{-u/B_j}}\right\}du=\int_{0}^{\infty}e^{-xu}uP(u)du.
\end{multline*}
The integral converges around zero due to asymptotic formula (\ref{eq:Puzeroasymp}) below. Convergence at infinity is obvious.  Hence, $(\log{W})''$ is c.m. iff $P(u)\ge0$ by Bernstein's theorem \cite[Theorem~1.4]{SSV}. $\hfill\square$

\rem By substituting $t=e^{-u}$ condition (\ref{eq:Pu-def}) can be
also written in the form
\begin{equation}\label{eq:Qu-def}
Q(t)=\sum\limits_{i=1}^{p}\frac{t^{a_i/A_i}}{1-t^{1/A_i}}-\sum\limits_{j=1}^{q}\frac{t^{b_j/B_j}}{1-t^{1/B_j}}\ge0~~\text{for}~t\in(0,1).
\end{equation}

If $(\log{W})''$ is completely monotonic two interesting possibilities appear:  $(\log{W})'>0$ and $(\log{W})'<0$.
These two possibilities are handled in the next two theorems.  Recall that a nonnegative function $f$ on $(0,\infty)$ is called Bernstein function if $f'$ is c.m. \cite[Definition~3.1 and~below]{SSV}.

\begin{theorem}\label{th:g-bernstein}
Suppose $\A,\B>0$ and $\a,\b\ge0$.  Then the function
$$
(\log{W})'=\sum\limits_{i=1}^{p}A_i\psi(A_ix+a_i)-\sum\limits_{j=1}^{q}B_j\psi(B_jx+b_j)
$$
is a Bernstein function iff condition \emph{(\ref{eq:Pu-def})}  holds and
$$
\lim_{x\to0}\left\{\sum\nolimits_{i=1}^{p}A_i\psi(a_i)-\sum\nolimits_{j=1}^{q}B_j\psi(b_j)\right\}\ge0.
$$
\end{theorem}
\textbf{Proof.} The proof is immediate from Lemma~\ref{lm:lcm-derivative} and the following equivalences
\begin{multline}\label{eq:cmequiv}
f \text{ is Bernstein function } \Leftrightarrow  f '\text{ is c.m. and }  f\ge0 \text{ on } (0,\infty)
\\
\Leftrightarrow f' \text{ is c.m. on } (0,\infty) \text{ and } \lim\limits_{x\to 0} f(x)\ge0.~~~\square
\end{multline}

\begin{theorem}\label{th:lcm-iff}
The function $W$ is l.c.m. if and only if
\begin{equation}\label{eq:nec1}
\sum\limits_{j=1}^{q}B_j=\sum\limits_{i=1}^{p}A_i,~~~\rho=\prod\limits_{i=1}^{p}A_i^{A_i}\prod\limits_{j=1}^{q}B_j^{-B_j}\le1
\end{equation}
and condition \emph{(\ref{eq:Pu-def})} holds. In the affirmative
case
\begin{equation}\label{eq:f-represent}
-(\log{W})'=\int_{0}^{\infty}e^{-xu}P(u)du+\log(1/\rho).
\end{equation}
\end{theorem}
\textbf{Proof.} We have
$$
-(\log{W})'=\sum\limits_{j=1}^{q}B_j\psi(B_jx+b_j)-\sum\limits_{i=1}^{p}A_i\psi(A_ix+a_i).
$$
Further it follows from (\ref{eq:psi-asimptotic}) that
$$
\psi(Cz+c)=\log{z} +\log{C}+\left(c-\frac{1}{2}\right)\frac{1}{Cz}+O(z^{-2})~~\text{as}~~z\to\infty.
$$
Substituting the latter formula into the former, we get
\begin{multline}\label{eq:LogWasymp}
-(\log{W})'=\log(x)\left\{\sum\limits_{j=1}^{q}B_j-\sum\limits_{i=1}^{p}A_i\right\}
+\left\{\sum\limits_{j=1}^{q}B_j\log(B_j)- \sum\limits_{i=1}^{p}A_i\log(A_i)\right\}\\
+\frac{1}{x}\left\{\sum\limits_{j=1}^{q}b_j-\sum\limits_{i=1}^{p}a_i+\frac{p-q}{2}\right\}
+O(x^{-2})~~\text{as}~~x\to\infty.
\end{multline}
If conditions (\ref{eq:Pu-def}) and (\ref{eq:nec1}) hold  than
$$
\lim\limits_{x\to\infty}(-\log W(x))'=-\log\rho\geq 0.
$$
Then Lemma~\ref{lm:lcm-derivative} and relations (\ref{eq:lcm}) imply that $W$ is l.c.m.
In the opposite direction, if  $W$ is l.c.m. than $-(\log{W})'$ is c.m. and hence decreasing which in view of (\ref{eq:LogWasymp}) implies that
$\sum_{i=1}^{p}{A_i}\ge\sum_{i=1}^{q}{B_i}$.  On the other hand, relations (\ref{eq:lcm}) show that $\lim_{x\to\infty}(-\log W(x))'\ge0$ so that $\sum_{i=1}^{p}{A_i}\le\sum_{i=1}^{q}{B_i}$ again by (\ref{eq:LogWasymp}). This proves the first condition in (\ref{eq:nec1}). Once this condition is satisfied nonnegativity of the limit leads to the second condition in (\ref{eq:nec1}). Finally, integrating (\ref{eq:represent 1}) we get
\begin{multline*}
(\log W(x))' =\int_{\infty}^x (\log W(x))'' dx+\log\rho
=\int\limits_{0}^{\infty}uP(u)du\int^x_{\infty}e^{-tu}dt+\log\rho
\\
=-\int\limits_{0}^{\infty}e^{-xu}P(u)du-\log(1/\rho).~~~~~\square
\end{multline*}

\rem Expression of the form $S_{\A}=\sum_{i=1}^{p}A_i\log(A_i)$ for positive numbers $A_i$ is known as Shannon's entropy in information theory, so that the second condition in (\ref{eq:nec1}) can be restated as ''entropy of the vector $\A$ does not exceed that of the vector $\B$''.

We now collect  the necessary conditions for logarithmic complete
monotonicity of $W$ in the next corollary.

\begin{corollary}\label{cr:necessary}
The following conditions are necessary for $W$ to be  logarithmically completely monotone:\\[1pt]

\emph{a)} $\sum_{j=1}^{q}B_j=\sum_{i=1}^{p}A_i$\\[1pt]

\emph{b)} $\rho\le1$\\[1pt]

\emph{c)} $\mu=\sum_{j=1}^{q}b_j-\sum_{i=1}^{p}a_i+\frac{1}{2}(p-q)\geq0$\\[1pt]

\emph{d)} $\min\limits_{1\leq{i\leq{p}}}(a_i/A_i)\leq\min\limits_{1\leq{j\leq{q}}}(b_j/B_j)$\\[1pt]
\end{corollary}

\textbf{Proof.}  Necessity of a) and b) has been demonstrated in Theorem~\ref{th:lcm-iff}.  Next, straightforward computation shows that
\begin{equation}\label{eq:Puzeroasymp}
P(u)=\frac{1}{u}\left\{\sum_{i=1}^pA_i-\sum_{j=1}^qB_j\right\}
+\left\{\sum\limits_{j=1}^{q}b_j-\sum\limits_{i=1}^{p}a_i+\frac{1}{2}(p-q)\right\}+O(u)~~\text{as}~~u\to0.
\end{equation}
Since the first sum is zero by a) we conclude that c) is necessary for nonnegativity of $P(u)$.  Further, using representation (\ref{eq:Qu-def}), we easily compute
$$
Q(t)=t^{\alpha}(1+o(1))-t^{\beta}(1+o(1))~~\text{as}~~t\to0,
$$
where $\alpha=\min\limits_{1\leq{i\leq{p}}}(a_i/A_i)$ and
$\beta=\min\limits_{1\leq{j\leq{q}}}(b_j/B_j)$.  Hence, we conclude that condition d) is also necessary.$\hfill\square$

\rem  Another method to show the necessity of condition $\mu\ge0$ (condition c) above) is as follows.  Calculating further  derivatives of $f$ and using
$$
(-1)^{n+1}\psi^{(n)}(x)=\frac{(n-1)!}{x^n}+\frac{n!}{2x^{n+1}}+O(x^{-n-2})~~\text{as}~~x\to\infty,
$$
we get
\begin{multline*}
(-1)^nf^{(n)}(x)=(-1)^{n+1}\sum\limits_{i=1}^{p}(A_i^{n+1}\psi^{(n)}(A_ix+a_i)
-(-1)^{n+1}\sum\limits_{j=1}^{q}B_j^{n+1}\psi^{(n)}(B_jx+b_j))
\\
=\sum\limits_{i=1}^{p}\left\{\frac{(n-1)!A_i^{n+1}}{(A_ix+a_i)^n}+\frac{n!A_i^{n+1}}{2(A_ix+a_i)^{n+1}}\right\}
-\sum\limits_{j=1}^{q}\left\{\frac{(n-1)!B_j^{n+1}}{(B_jx+b_j)^n}+\frac{n!B_j^{n+1}}{2(B_jx+b_j)^{n+1}}\right\}+O(x^{-n-2})
\\
=\sum\limits_{i=1}^{p}\left\{\frac{(n-1)!A_i}{x^n}\left(1-\frac{na_i}{A_ix}\right)+\frac{n!}{2x^{n+1}}\right\}
-\sum\limits_{j=1}^{q}\left\{\frac{(n-1)!B_j}{x^{n}}\left(1-\frac{nb_j}{B_{j}x}\right)+\frac{n!}{2x^{n+1}}\right\}+O(x^{-n-2})
\\
=\frac{(n-1)!}{x^n}\left\{\sum\limits_{i=1}^{p}A_i-\sum\limits_{j=1}^{q}B_j\right\}
+\frac{n!}{x^{n+1}}\left\{\sum\limits_{j=1}^{q}b_j-\sum\limits_{i=1}^{p}a_i+\frac{1}{2}(p-q)\right\}+O(x^{-n-2})
\\
=\frac{n!}{x^{n+1}}\left\{\sum\limits_{j=1}^{q}b_j-\sum\limits_{i=1}^{p}a_i+\frac{1}{2}(p-q)\right\}+O(x^{-n-2}),
\end{multline*}
where the first sum vanishes due to the necessary condition a) from Corollary~\ref{cr:necessary}.  This shows that this must require that $\sum_{j=1}^{q}b_j-\sum_{i=1}^{p}a_i+\frac{1}{2}(p-q)\geq0$ in order for $(-1)^nf^{(n)}(x)$ to be nonnegative for large $x$.

Given the vectors $\A,\B,\a,\b$ condition (\ref{eq:Pu-def}) is not easy to verify. The next theorem provides some practical
sufficient conditions.
\begin{theorem}\label{th:Pu-sufficient}
Inequality \emph{(\ref{eq:Pu-def})} is true  if any of the following conditions  holds\emph{:}

\bigskip

\noindent\emph{(a)} $\sum_{j=1}^{q}B_j=\sum_{i=1}^{p}A_i$ and
$\max\limits_{1\le{i}\le{p}}(a_i/A_i)\leq \min\limits_{1\le{j}\le{q}}(b_j-1)/B_j$\emph{;}

\bigskip

\noindent\emph{(b)} $p=q$,
$
\sum_{i=1}^{p}B_i=\sum_{i=1}^{p}A_i~\text{with}~A_i\ge{B_i}~\text{for}~i=1,\ldots,p-1,
$
and
$$
\max\limits_{1\le{k}\le p-1}b_k/B_k\le(b_{p}-1)/B_p,~~~a_i/A_i\le(b_i-1)/B_i
~~\text{for}~i=1,\ldots,p;
$$

\smallskip

\noindent\emph{(c)} $p=q$,  $0\leq {a_1}/{A_1}\leq {a_2}/{A_2}\leq\ldots\leq {a_p}/{A_p}$,
$0\leq {b_1}/{B_1}\leq {b_2}/{B_2}\leq\ldots \leq{b_p}/{B_p}$,\\[7pt]
\indent $0<{1}/{A_1}\leq {1}/{A_2}\leq\ldots\leq {1}/{A_p}$, $0<{1}/{B_1}\leq {1}/{B_2}\leq\ldots\leq {1}/{B_p}$
and\\[7pt]
$\sum_{i=1}^{k}a_i/A_i\leq \sum_{i=1}^{k}b_i/B_i$,
$\sum_{i=1}^{k}1/A_i\leq\sum_{i=1}^{k}1/B_i$ for each integer $k=1,\ldots,p$.
\end{theorem}

\textbf{Proof.} (a) Denote $a'_i=a_i/A_i$, $A'_i=1/A_i,$ $b'_i=b_i/B_i$, $B'_i=1/B_i,$ $x=e^u.$ Then by the mean value theorem we have
$$
\sum\limits_{i=1}^{p}\frac{e^{-a_iu/A_i}}{1-e^{-u/A_i}}=\sum\limits_{i=1}^{p}\frac{e^{-a_i'u}}{1-e^{-A_i'u}}=
\sum\limits_{i=1}^{p}\frac{1}{x^{a_i'}-x^{a_i'-A'_i}}=\sum\limits_{i=1}^{p}
\frac{1}{A'_ix^{c_i}\log{x}}=\frac{1}{\log{x}}\sum\limits_{i=1}^{p}A_ix^{-c_i},
$$
where $c_i\in(a_i'-A_i', a'_i)$. Hence, for $P(u)$ we get
\begin{equation}\label{eq:factorp}
P(u)=\frac{1}{u}\left\{\sum\limits_{i=1}^{p}A_ie^{-c_iu}
-\sum\limits_{j=1}^{q}B_je^{-d_ju}\right\},
\end{equation}
where $d_j$ is some point of the interval $(b_j'-B_j', B'_j)$.
Denote $c=\max{c_i}$, $d=\min{d_j}$. Suppose $\max(a_i/A_i)\leq\min(b_j-1)/B_j$, which is equivalent to
$\max(a_i')\leq \min(b_j'-B_j')$. Then $c<d$ and nonnegativity of $P(u)$ follows from the relations
$$
\sum\limits_{i=1}^{p} A_ix^{-c_i}\geq\sum\limits_{i=1}^{p} A_ix^{-c}\geq
\sum\limits_{i=1}^{p} A_ix^{-d}=\sum\limits_{j=1}^{q}
B_jx^{-d}\geq\sum\limits_{j=1}^{q} B_jx^{-d_j}.
$$

(b) Denote $S=\sum\limits_{j=1}^{p}B_j=\sum\limits_{i=1}^{p}A_i$.  Nonnegativity of $P(u)$ follows from representation
(\ref{eq:factorp}) and the following chain
\begin{multline*}
\sum\limits_{i=1}^{p}A_ie^{-c_iu}\geq\sum\limits_{i=1}^{p}A_ie^{-d_iu}
=\sum\limits_{i=1}^{p-1}A_ie^{-d_iu}+\left(\!S-\sum\limits_{i=1}^{p-1} A_i\!\right)e^{-d_pu}
\\
=\sum\limits_{i=1}^{p-1} A_i(e^{-d_iu}-e^{-d_p u})+Se^{-d_pu}
\ge\sum\limits_{i=1}^{p-1} B_i(e^{-d_iu}-e^{-d_p u})+Se^{-d_pu}=\sum\limits_{i=1}^{p} B_ie^{-d_i u}.
\end{multline*}
The first inequality here is due to conditions
$a_i/A_i\le(b_i-1)/B_i$, $i=1,\ldots,p$ which imply that
$c_i<d_i$. The second inequality follows from conditions
$A_i\ge{B_i}$, $i=1,\ldots,p-1$, combined with inequality
$\max\limits_{1\le{k}\le{p-1}}b_k/B_k\le(b_{p}-1)/B_p$ which
ensures that each term $e^{-d_iu}-e^{-d_p u}$ is nonnegative.

(c) Expanding each expression like $(1-e^{-u/A_i})^{-1}$ in geometric series and exchanging the order of summations we get
$$
P(u)=\sum\limits_{l=0}^{\infty}\sum\limits_{i=1}^{p}\left(e^{-a_iu/A_i-lu/A_i}-e^{-b_iu/B_i-lu/B_i}\right).
$$
Conditions c) guarantee that for any $l\ge0$, $u\ge0$ we have
$$
(b_1 u/B_1+lu/B_1,...,b_p u/B_p+lu/B_p)\prec^{W}(a_1u/A_1+lu/A_1,...,a_p u/A_p+lu/A_p)
$$
(see (\ref{eq:amajorb}) for definition of weak supermajorization
$\prec^{W}$).  Hence, according to \cite[Proposition~4.B.2]{MOA}
we conclude that
$$
\sum\limits_{i=1}^{p}\left(e^{-a_iu/A_i-lu/A_i}-e^{-b_iu/B_i-lu/B_i}\right)\geq0
$$
since $x\to{e^{-x}}$ is convex and decreasing.$\hfill\square$

\rem  Conditions (c) from Theorem~\ref{th:Pu-sufficient} are only compatible with condition $\sum_{j=1}^{p}A_j=\sum_{j=1}^{p}B_j$ if $\A=\B$ (up to permutation).  Indeed, majorization $\B'\prec^W\A\!'$ (recall that $A_j'=A_j^{-1}$, $B_j'=B_j^{-1}$) forms a part of conditions (c). According to \cite[Proposition~4.B.2]{MOA} this majorization implies that
$$
\sum\limits_{j=1}^{p}B_j=\sum\limits_{j=1}^{p}\frac{1}{B'_j}\leq\sum\limits_{j=1}^{p}\frac{1}{A'_j}=\sum\limits_{j=1}^{p}A_j,
$$
because the function $x\to1/x$ is decreasing and convex. Further, \cite[3.A.6.a]{MOA} says that for continuous strictly decreasing functions the above inequality is strict unless $\A\!'$ is a permutation of $\B'$. This brings us to the conclusion that for $\A\neq\B$ (modulo permutations) conditions (c) can be used to check whether $(\log{W})'$ is a Bernstein function using Theorem~\ref{th:g-bernstein}, but they \emph{cannot} be used to check whether $W$ is completely monotone. On the other hand, if $\A=\B$, conditions (c) reduce to checking the majorization $\a'\prec^W\b'$ ($a_j'=a_jA_j^{-1}$, $b_j'=b_jA_j^{-1}$). This majorization is weaker than (a) or (b) in this particular situation.  One reason why the case $\A=\B$ might be important is applications in probability as explained  in Remark~\ref{rm:randombeta} below.
\bigskip

\textbf{Example~2.} The function
 $$
 x\to\frac{\Gamma(2x+0.4)\Gamma(3x+2.4)\Gamma(x+0.9)}{\Gamma(x+2)\Gamma(5x+6)}
 $$
 is logarithmically completely monotone since $A_1+A_2+A_3=B_1+B_2$, $\rho=0.03456$  and $\max_{1\le{i}\le3}(a_i/A_i)=0.9$, $\min_{1\le{i}\le2}((b_i-1)/B_i)=1$. Hence, necessary conditions (\ref{eq:nec1}) are satisfied and $P(u)\ge0$ according to Theorem~\ref{th:Pu-sufficient}(a).  Note, that  neither conditions (b) nor (c) from Theorem~\ref{th:Pu-sufficient} can be applied here as $p\neq{q}$.

\textbf{Example~3.} The function
$$
x\to\frac{\Gamma(2x+0.8)\Gamma(3x+8)\Gamma(1.4x+2.3)}{\Gamma(x+1.5)\Gamma(2.4x+7.8)\Gamma(3x+11)}
$$
is logarithmically completely monotone since $A_1+A_2+A_3=B_1+B_2+B_3$, $\rho=0.783668$, so that necessary conditions (\ref{eq:nec1}) are satisfied and $A_1>B_1$, $A_2>B_2$,
$\max_{1\le{i}\le2}b_i/B_i=3.25<(b_3-1)/B_3=3.33$, $a_1/A_1=0.4<(b_1-1)/B_1=0.5$, $a_2/A_2=2.66<(b_2-1)/B_2=2.83$,
$a_3/A_3=1.643<(b_3-1)/B_3=3.33$, so that $P(u)\ge0$ by Theorem~\ref{th:Pu-sufficient}(b).  Note that conditions (a) and (c) of Theorem~\ref{th:Pu-sufficient} are violated here.

\textbf{Example~4.} The function
$$
x\to 4\psi(4x+0.7)+2\psi(2x+1.8)-3\psi(3x+0.6)-\psi(x+1.2)
$$
is a Bernstein function, since $A_1=4>A_2=2$, $B_1=3>B_2=1$, $A_1^{-1}=0.25<B_1^{-1}=0.33$,
$A_1^{-1}+A_2^{-1}=0.75<B_1^{-1}+B_2^{-1}=1.33$, $a_1A_1^{-1}=0.175<b_1B_1^{-1}=0.2$,
$a_1A_1^{-1}+a_2A_2^{-1}=1.075<b_1B_1^{-1}+b_2B_2^{-1}=1.4$, so that $P(u)\ge0$ by Theorem~\ref{th:Pu-sufficient}(c) and the claim follows from Corollary~\ref{cr:g-bernstein}.  Note, that the function $W$ with parameters from this example is not logarithmically completely monotone since necessary conditions (\ref{eq:nec1}) are violated.

\textbf{Example~5.} The function
$$
x\to\frac{\Gamma(3x+0.8)\Gamma(2.2x+1.8)\Gamma(1.4x+2.3)}{\Gamma(3x+1.2)\Gamma(2.2x+1.7)\Gamma(1.4x+2.5)}
$$
is logarithmically completely monotone by Theorem~\ref{th:Pu-sufficient}(c).   Indeed,  $A_i=B_i$, $i=1,2,3$, and majorization conditions are easily  verified.
\bigskip

\paragraph{4. The representing measure.}
We will need a particular case of Fox's $H$-function defined by
\begin{equation}\label{eq:Fox}
H_{q,p}^{p,0}\left(z\left|\begin{array}{l} (\B,\b)\\
(\A,\a)\end{array}\right.\right)=
\frac{1}{2\pi{i}}\int\limits_{\L}
\frac{\prod\nolimits_{k=1}^{p}\Gamma(A_ks+a_k)}
{\prod\nolimits_{j=1}^{q}\Gamma(B_j s+b_j)} z^{-s}ds,
\end{equation}
where $A_k,B_j>0$ and $a_k,b_j$ are real.  The contour $\L$ can be either the left loop $\L_{-}$ starting at $-\infty+i\alpha$ and ending at $-\infty+i\beta$ for some $\alpha<0<\beta$ such that all poles of the integrand lie inside the loop; or the right loop $\L_{+}$ starting at $\infty+i\alpha$ and ending at  $\infty+i\beta$ and leaving all poles on the left; or the vertical line $\L_{ic}$, $\Re{z}=c$,  traversed upwards and leaving all poles of the integrand on the left.  Denote the rightmost pole of the integrand by $\gamma$:
$$
\gamma=-\min\limits_{1\leq k\leq p}({a_k}/{A_k}).
$$
Recall the definition of $\rho$ from (\ref{eq:nec1}) and the definition of $\mu$ from Corollary~1(c):
$$
\rho=\prod\limits_{k=1}^{p}A_k^{A_k}\prod\limits_{j=1}^{q}B_j^{-B_j},
~~~\mu=\sum\limits_{j=1}^{q}b_j-\sum\limits_{k=1}^{p}a_k+\frac{p-q}{2}.
$$
Existence conditions of $H$-function under each choice of the
contour $\L$ have been thoroughly considered in the book
\cite{KilSaig}. Under restrictions (\ref{eq:nec1}) and $x>0$
Theorem~1.1 from \cite{KilSaig} states that $H_{q,p}^{p,0}(x)$
exists if we choose $\L=\L_{+}$ or $\L=\L_{ic}$ under additional
restriction $\mu>1$. Only the second choice of the contour ensures
the existence of the Mellin transform of  $H_{q,p}^{p,0}(x)$ as
demonstrated in \cite[Theorem~2.2]{KilSaig}. In our next theorem
we relax the condition $\mu>1$ to $\mu>0$ and demonstrate that the
first condition in (\ref{eq:nec1}) leads to the compact support of
$H_{q,p}^{p,0}(x)$.
\begin{theorem}\label{th:Fox}
Suppose $\mu>0$ and $\sum_{i=1}^pA_i=\sum_{j=1}^q B_j$.  Then  the integral over the contour $\L_{ic}$ with $c>\gamma$ in the definition of $H_{q,p}^{p,0}(x)$ converges for $x>0$ except possibly for $x=\rho$ and
\begin{equation}\label{eq:H-vanish}
H_{q,p}^{p,0}\left(x\left|\begin{array}{l}(\B,\b)\\(\A,\a)\end{array}\right.\right)=0
\end{equation}
for $x>\rho$. Moreover, under these restrictions, the Mellin transform exits for $\Re{s}>\gamma$ and
$$
\int_{0}^{\rho}H_{q,p}^{p,0}\left(x\left|\begin{array}{l}(\B,\b)\\(\A,\a)\end{array}\right.\right)x^{s-1}dx
=\frac{\prod\nolimits_{k=1}^{p}\Gamma(A_ks+a_k)}{\prod\nolimits_{j=1}^{q}\Gamma(B_js+b_j)}.
$$
\end{theorem}
\textbf{Proof.}  By applying Stirling's asymptotic formula and in
view of $\sum_{k=1}^p A_k=\sum_{j=1}^q B_j$ it is not difficult to
derive the formula \cite[(2.2.4)]{ParKam}
\begin{equation}\label{eq:gamm}
\frac{\prod\nolimits_{k=1}^{p}\Gamma(A_ks+a_k)}{\prod\nolimits_{j=1}^{q}\Gamma(B_js+b_j)}=
A\rho^{s}s^{-\mu}+\rho^sg(s),
\end{equation}
with $g(s)=O(s^{-\mu-1})$ as $|s|\to\infty$ uniformly in $|\arg{s}|\le\pi-\epsilon$ for any $0<\epsilon<\pi$ and
$$
A=(2\pi)^{(p-q)/2}\prod\nolimits_{k=1}^{p}A_k^{a_k-1/2}\prod\nolimits_{j=1}^{q}B_k^{1/2-b_j}.
$$
This asymptotic behavior implies that $t\to{g(c+it)}$ is absolutely integrable continuous function on the real line for any $c>\gamma$ so that the integral
$v(x)=\frac{1}{2\pi{i}}\int_{c-i\infty}^{c+i\infty}(x/\rho)^{-s}g(s)ds$ exists and we are in the position to apply the Mellin inversion theorem yielding
$$
\int_{0}^{\infty}x^{s-1}v(x)dx=\rho^sg(s).
$$
Denote
$$
h(x)= \frac{1}{2\pi{i}}\int\limits_{c-i\infty}^{c+i\infty}A\rho^{s}s^{-\mu}x^{-s}ds.
$$
An application of \cite[(6), \S12]{Bochner} after a simple rearrangement gives
$$
h(x)=
\frac{A}{2\pi}\int\limits_{-\infty}^{+\infty}\frac{e^{\log(\rho/x)(c+it)}}{(c+it)^\mu}dt=
\left\{\!\!\!\begin{array}{l}\dfrac{A}{\Gamma(\mu)}\left(\log\frac{\rho}{x}\right)^{\mu-1}\!\!\!,~0<x<\rho,\\[10pt]
0,\ x>\rho.\end{array}\right.
$$
 It then follows from (\ref{eq:gamm}) that
$$
H_{q,p}^{p,0}\left(x\left|\begin{array}{l}
(\B,\b)\\(\A,\a)\end{array}\right.\right)= \frac{1}{2\pi{i}}\int\limits_{c-i\infty}^{c+i\infty}
\frac{\prod\nolimits_{k=1}^{p}\Gamma(A_ks+a_k)}{\prod\nolimits_{j=1}^{q}\Gamma(B_js+b_j)} x^{-s}ds=h(x)+v(x)
$$
and the above integral exists for all $x>0$, $x\neq\rho$.  Direct computation gives
$$
\int\limits_{0}^{\infty}x^{s-1}H_{q,p}^{p,0}\left(x\left|\begin{array}{l}
(\B,\b)\\(\A,\a)\end{array}\right.\right)dx=\frac{A}{\Gamma(\mu)}\int\limits_{0}^{\rho}
\left(\log\frac{\rho}{x}\right)^{\mu-1}x^{s-1}dx+\rho^sg(s)=A\rho^{s}s^{-\mu}+\rho^sg(s),
$$
which is equivalent to
$$
\int\limits_{0}^{\infty}x^{s-1}H_{q,p}^{p,0}\left(x\left|\begin{array}{l}
(\B,\b)\\(\A,\a)\end{array}\right.\right)dx
=\frac{\prod\nolimits_{k=1}^{p}\Gamma(A_ks+a_k)}{\prod\nolimits_{j=1}^{q}\Gamma(B_js+b_j)}.
$$
Next we show that (\ref{eq:H-vanish}) holds so that the integration in the above formula is in fact over the interval $(0,\rho)$. To this end we apply an argument similar to Jordan's lemma. According to the definition (\ref{eq:Fox})
$$
H_{q,p}^{p,0}\left(z\left|\begin{array}{l}
(\B,\b)\\(\A,\a)\end{array}\right.\right)=
\lim\limits_{R\to\infty}\frac{1}{2\pi i}\int\limits_{c-iR}^{c+iR}
\frac{\prod\nolimits_{k=1}^{p}\Gamma(A_ks+a_k)}{\prod\nolimits_{j=1}^{q}\Gamma(B_j
s+b_j)} z^{-s}ds.
$$
By Cauchy's theorem the last integral equals to the integral along the right semicircle of radius $R$ centered at $c$.
Hence, we need to prove that
$$
I_R=\frac{R}{2\pi }\int\limits_{-\pi/2}^{\pi/2}
\frac{\prod\nolimits_{k=1}^{p}\Gamma(A_kc+A_k Re^{i\varphi}+a_k)}
{\prod\nolimits_{j=1}^{q}\Gamma(B_j c+B_j Re^{i\varphi}+b_j)}
e^{-(c+R\cos\varphi+iR\sin\varphi)\log x+i\varphi}d\varphi
$$
goes to zero as  $R\to\infty$ for $x>\rho$. Setting $s=Re^{i\varphi}$ in
(\ref{eq:gamm}) we get
$$
\frac{\prod\nolimits_{k=1}^{p}\Gamma(A_kc+A_kRe^{i\varphi}+a_k)}
{\prod\nolimits_{j=1}^{q}\Gamma(B_j c+B_jRe^{i\varphi}+b_j)}
=A\rho^{Re^{i\varphi}}(Re^{i\varphi})^{-\mu}(1+O(R^{-1}))~~\text{as}~~R\to\infty.
$$
Hence, for sufficiently large $R$ and some constant $M_1$
$$
|I_R|\leq \frac{M_1e^{-c\log
x}}{2\pi}\int\limits_{-\pi/2}^{\pi/2}R^{-\mu+1}e^{-R\log\frac{x}{\rho}\cos\varphi}d\varphi.
$$
Employing the elementary inequality
$\cos\varphi\geq1-\frac{2}{\pi}\varphi$, $0<\varphi<\pi/2$, we
obtain the following estimate with some constant $M_2$ (recall that $x>\rho$):
$$
|I_R|\leq2M_2R^{-\mu+1}\int\limits_{0}^{\pi/2}e^{-R(1-\frac{2}{\pi}\varphi)\log\frac{x}{\rho}}d\varphi
=\pi{M_2}R^{-\mu}\frac{(1-e^{-R\log\frac{x}{\rho}})}{\log\frac{x}{\rho}}.
$$
Hence, $\lim\limits_{R\to\infty}I_R=0$ which completes the proof of the theorem.$\hfill\square$

\begin{theorem}\label{th:measure}
Suppose that $\mu>0$ and
$$
W(x)=\frac{\prod\nolimits_{i=1}^{p}\Gamma(A_ix+a_i)}{\prod\nolimits_{j=1}^{q}\Gamma(B_jx+b_j)}
$$
is l.c.m. Then
$$
W(x)=\int_{\log(1/\rho)}^{\infty}e^{-tx}H_{q,p}^{p,0}\left(e^{-t}\left|\begin{array}{l}(\B,\b)\\(\A,\a)\end{array}\right.\right)\!dt,
$$
and the $H$-function in the integrand is nonnegative. In particular, the conclusion is true if conditions \emph{(\ref{eq:nec1})} and hypotheses of Theorem~\emph{\ref{th:Pu-sufficient}} are satisfied.
\end{theorem}

\textbf{Proof.} According to Theorem~\ref{th:lcm-iff} logarithmic complete monotonicity of $W(x)$ implies
$\sum_{k=1}^{p}A_k=\sum_{j=1}^{q}B_j$, so that we are in the position to apply Theorem~\ref{th:Fox} yielding the formula
$$
W(x)=\frac{\prod\nolimits_{i=1}^{p}\Gamma(A_ix+a_i)}{\prod\nolimits_{j=1}^{q}\Gamma(B_jx+b_j)}
=\int_{0}^{\rho}H_{q,p}^{p,0}\left(u\left|\begin{array}{l}(\B,\b)\\(\A,\a)\end{array}\right.\right)u^{x-1}du.
$$
The claimed Laplace transform representation for $W(x)$ follows by substitution $u=e^{-t}$.
Nonnegativity of $H_{q,p}^{p,0}(u)$ on $(0,\rho)$ follows from Bernstein theorem in view of uniqueness of the measure with given Laplace transform, see \cite[Proposition~1.2]{SSV} or  \cite[Theorem~6.3]{Widder}. $\hfill\square$

Theorem~\ref{th:measure} requires $\mu$ to be strictly positive, while the necessary conditions from Corollary~\ref{cr:necessary} allow for $\mu=0$.  In analogy with (\ref{eq:psi0measure}) in that case we believe in the validity of the following conjecture.

\begin{conjecture}\label{cj:mu0measure}
For $\mu=0$ the representing measure is given by
$$
W(x)=\int_{[\log(1/\rho),\infty)}e^{-tx}\left\{A\delta_{\log(1/\rho)}+H_{q,p}^{p,0}\left(e^{-t}\left|\begin{array}{l}(\B,\b)\\(\A,\a)\end{array}\right.\right)\!\right\}dt,
$$
where $A$ is a positive constant and $\delta_{\log(1/\rho)}$ denotes the unit mass concentrated at the point $\log(1/\rho)$.
\end{conjecture}

Unlike the case of $G$-function (\ref{eq:psi0measure}), no expansion similar to N{\o}rlund's  is available for $H$-function.
We plan to study its behavior in the neighborhood of the singular point $\log(1/\rho)$ in our future work.

\begin{corollary}\label{cr:H-infdiv}
Suppose $\mu>0$, $\a,\b>0$ and conditions \emph{(\ref{eq:Pu-def})} and \emph{(\ref{eq:nec1})} are satisfied \emph{(}in particular, \emph{(\ref{eq:Pu-def})} is true under Theorem~\emph{\ref{th:Pu-sufficient}}\emph{)}. Then the function
\begin{equation}\label{eq:H-infdiv}
t\to \frac{\prod_{j=1}^{q}\Gamma(b_j)}{\prod_{i=1}^{p}\Gamma(a_i)}H_{q,p}^{p,0}\left(e^{-t}\left|\begin{array}{l}(\B,\b)\\(\A,\a)\end{array}\right.\right)
\end{equation}
represents an infinitely divisible probability density on $(0,\infty)$.
\end{corollary}
\textbf{Proof.}  The result follows from Theorem~\ref{th:measure} combined with \cite[Definition~5.6]{SSV} and \cite[Theorem~5.9]{SSV}.  See also \cite[Theorem~7]{Koumandos}.  $\hfill\square$

\rem\label{rm:randombeta} The case $\A=\B$ is also connected to probability as follows. Suppose $\zeta_k$, $k=1,\ldots,p$ are independent beta-distributed random variables, so that  $\zeta_k$ has the density
$$
\frac{\Gamma(\alpha_k+\beta_k)}{\Gamma(\alpha_k)\Gamma(\beta_k)}t^{\alpha_k}(1-t)^{\beta_k},~~t\in(0,1).
$$
Then the random variable $u=\prod_{k=1}^{p}\zeta_k^{A_k}$ with $A_k>0$, has the following moments
$$
E(u^{x-1})=\prod\limits_{k=1}^{p}\frac{\Gamma(\alpha_k+\beta_k)}{\Gamma(\alpha_k)}
\prod\limits_{k=1}^{p}\frac{\Gamma(A_kx+\alpha_k-A_k)}{\Gamma(A_kx+\alpha_k+\beta_k-A_k)},~~~x>0.
$$
The probability density of  $u$ is expressed via Fox's $H$-function:
$$
\prod_{k=1}^{p}\frac{\Gamma(\alpha_k+\beta_k)}{\Gamma(\alpha_k)}H_{p,p}^{p,0}\left(z\left|\begin{array}{l} (A_j,\alpha_j+\beta_j-A_j),j=1,\ldots,p\\
(A_j,\alpha_j-A_j),j=1,\ldots,p\end{array}\right.\right)
$$
See \cite[Section~4.2.1]{MSH}.

\begin{theorem}\label{th:hrepresentation}
Suppose $\A,\B>0$, $\a,\b\ge0$, $\mu>0$, $\rho\le1$ and $\sum_{i=1}^{p}A_i=\sum_{j=1}^{q}B_j$.
Then for all $x\in(0,\rho)$ the following identity holds:
\begin{equation}\label{eq:H-identity}
H_{q,p}^{p,0}\left(x\left|\begin{array}{l}(\B,\b)\\(\A,\a)\end{array}\right.\right)=\frac{1}{\log(\rho/x)}\int_{x/\rho}^{1}
H_{q,p}^{p,0}\left(\frac{x}{u}\left|\begin{array}{l}(\B,\b)\\(\A,\a)\end{array}\right.\right)\frac{Q(u)}{u}du,
\end{equation}
where $Q(u)$ is defined in \emph{(\ref{eq:Qu-def})}.
\end{theorem}
\textbf{Proof.} Since $-\log(\rho)=S_{\B}-S_{\A}$, where $S_{\A}$ ($S_{\B}$) stands for the entropy of $\A$ ($\B$), we adopt the notation $\Delta{S}=S_{\B}-S_{\A}=-\log(\rho)$.  Further, denote
$$
I(t)=\left\{ \begin{array}{ll}
\!\!1, &t\geq \Delta{S}\\
\!\!0, &t< \Delta{S}
\end{array},\right.\
$$
and define the (signed) measure
$$
d\nu(u)=\Delta{S}\delta_0+P(u)du,
$$
where $\delta_0$ is a unit mass at zero.  It follows from the proof of Theorem~\ref{th:lcm-iff} and the asymptotic formula (\ref{eq:Puzeroasymp}) that representation ({\ref{eq:f-represent}}) is true under conditions of the theorem. In terms of the measure $d\nu$ ({\ref{eq:f-represent}}) takes the form:
$$
W'(s)=W(\log{W})'=-W\int_{[0,\infty)}e^{-su}d\nu(u).
$$
On the other hand,  under the hypotheses of Theorem~\ref{th:Fox} we have the representation
$$
W(s)=\int_{0}^{\infty}\!e^{-ts}H(e^{-t})I(t)dt,~\text{where}~H(x)=H_{q,p}^{p,0}\left(x\left|\begin{array}{l}(\B,\b)\\(\A,\a)\end{array}\right.\right).
$$
Hence,
$$
W'(s)=-\int_{0}^{\infty}\!te^{-ts}H(e^{-t})I(t)dt.
$$
Putting these representations together yields
$$
\int_{0}^{\infty}\!te^{-ts}H(e^{-t})I(t)dt=\int_{0}^{\infty}\!e^{-ts}H(e^{-t})I(t)dt\int_{[0,\infty)}e^{-su}d\nu(u).
$$
The convolution theorem for the Laplace transform \cite[Theorem~11.4]{Widder} then leads to the formula
$$
tH(e^{-t})I(t)=\int_{[0,t)}\!\!H(e^{-(t-\tau)})I(t-\tau)d\nu(\tau)~~~\text{for}~t\in(\Delta{S},\infty).
$$
Recalling the definition of $d\nu$ we can rewrite the above relation as
$$
tH(e^{-t})=\Delta{S}H(e^{-t})+\int_{0}^{t-\Delta{s}}H(e^{-t+\tau})P(\tau)d\tau
$$
or
$$
H(e^{-t})=\frac{1}{t-\Delta{S}}\int_{0}^{t-\Delta{S}}H(e^{-t+\tau})P(\tau)d\tau.
$$
Substituting $x=e^{-t}$ and $u=e^{-\tau}$ yields identity (\ref{eq:H-identity}).$\hfill\square$

Theorem~\ref{th:lcm-iff} shows that under additional restrictions (\ref{eq:nec1}) we have the implication
$$
Q(t)\ge0~\text{on}~[0,1)~\Rightarrow~H_{q,p}^{p,0}\left(x\left|\begin{array}{l}(\B,\b)\\(\A,\a)\end{array}\right.\right)\ge0~\text{on}~(0,\rho).
$$
The following much stronger assertion is supported by numerical evidence:

\begin{conjecture}\label{cj:zerosH}
Suppose $\A,\B>0$, $\a,\b\ge0$, $\mu\ge0$, $\rho\le1$ and $\sum_{i=1}^{p}A_i=\sum_{j=1}^{q}B_j$. Then
$$
\#\biggl\{\text{zeros of}~H_{q,p}^{p,0}\left(x\left|\begin{array}{l}(\B,\b)\\(\A,\a)\end{array}\right.\right)~\text{on}~(0,\rho)\biggr\}
\le \#\biggl\{\text{zeros of}~Q(t)~\text{on}~(0,1)\biggr\}.
$$
\end{conjecture}
Using Theorem~\ref{th:hrepresentation} we were able to demonstrate that this conjecture is true for the case when the left hand side is equal to one.

\bigskip

\textbf{5. Acknowledgements.} We thank the anonymous referee for correcting a mistake in original version of the paper and numerous useful suggestions that helped to improve the exposition substantially.  This work has been supported by the Russian Science Foundation under project 14-11-00022.

\end{document}